\documentclass[12pt]{article}                
\usepackage{graphicx}
\usepackage{psfrag}
\usepackage{amsmath, amssymb}
\usepackage{mathptmx}

\newcommand{\A}{{\cal A}}

\newcommand{\LL}{{\cal L}}

\newcommand{\p}{\mathbb{P}}
\renewcommand{\t}[1]{\widetilde{#1}}
\newcommand{\Reals}{\mathbb{R}}
\newcommand{\Natural}{\mathbb{N}}

\newcommand\qed{\hfill\hbox{\rlap{$\sqcap$}$\sqcup$}}

\newtheorem{proposition}{Proposition}
\newtheorem{lemma}{Lemma}
\newtheorem{theorem}{Theorem}

\begin{document}

\title{A hierarchical version of the de Finetti and Aldous-Hoover representations.}
\author{Tim Austin\thanks{Courant Institute, New York University, tim@cims.nyu.edu. Partially supported by fellowship from Clay Mathematics Institute} \and Dmitry Panchenko\thanks{Department of Mathematics, Texas A\&M University, panchenk@math.tamu.edu. Partially supported by NSF grant.}\\
}
\date{}
\maketitle
\begin{abstract}
We consider random arrays indexed by the leaves of an infinitary rooted tree of finite depth, with the distribution invariant under the rearrangements that preserve the tree structure. We call such arrays hierarchically exchangeable and prove that they satisfy an analogue of de Finetti's theorem. We also prove a more general result for arrays indexed by several trees, which includes a hierarchical version of the Aldous-Hoover representation.
\end{abstract} 
\vspace{0.5cm}
Key words: exchangeability, spin glasses.\\
Mathematics Subject Classification (2010): 60G09, 60K35

\section{Introduction}

The subject of exchangeability is prevalent in probability theory (see e.g. \cite{Aldous2}, Chapters $7$--$9$ in \cite{Kallenberg}, or \cite{Aldous3}, \cite{Austin} and \cite{Austin2} for recent overviews and results) and the goal of this paper is to study another notion of exchangeability that is motivated by spin glass models and, in particular, by the work of M\'ezard and Parisi on diluted models, \cite{Mezard}.

\begin{figure}[t]
 \centering
 \psfrag{a_0}{ {\small $\emptyset$}}\psfrag{a_1}{ {\small $n_1$}}\psfrag{a_2}{ {\small $(n_1,\ldots,n_{r-1})$}}\psfrag{AcodeR}{ {\small $\alpha = (n_1,\ldots,n_{r})$}}\psfrag{n_code1}{ {\small $\Natural^{1}$}}\psfrag{n_r1}{ {\small $\Natural^{r-1}$}}\psfrag{n_r}{ {\small $\Natural^r$}}
 \includegraphics[width=0.84\textwidth]{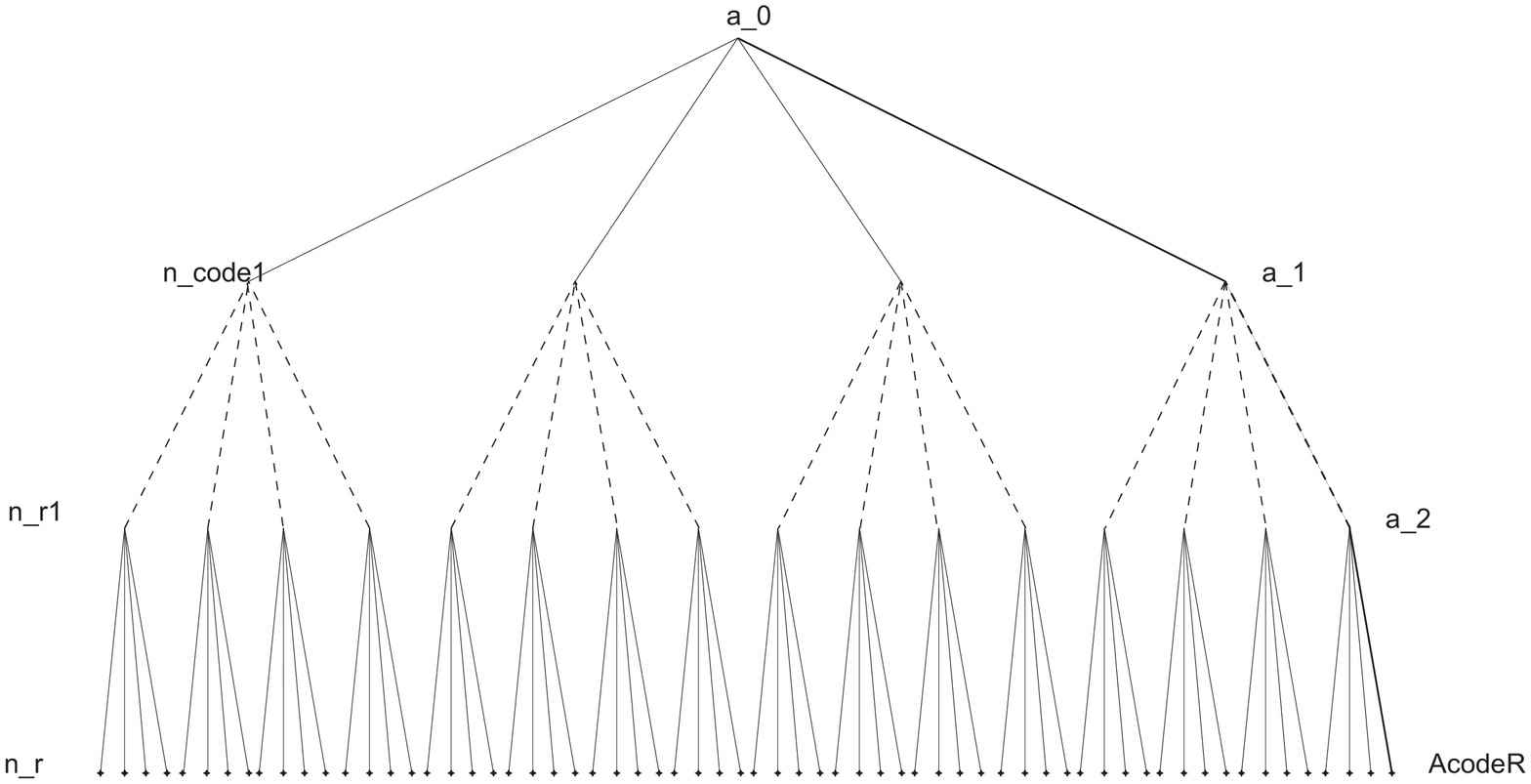}
 \caption{\label{Fig1} Index set $\Natural^r$ as the leaves of the infinitary tree $\A(r)$.}
 \end{figure}

We begin by considering an array $(X_\alpha)_{\alpha\in \Natural^r}$ of random variables $X_\alpha$ indexed by $\alpha\in \Natural^r$ for some integer $r\geq 1$, whose distribution is invariant under certain rearrangements of the indices. We will think of $\Natural^r$ as the set of leaves of a rooted tree (see Fig. \ref{Fig1}) with the vertex set
\begin{equation}
\A(r) = \Natural^0 \cup \Natural \cup \Natural^2 \cup \ldots \cup \Natural^r,
\label{Atree}
\end{equation}
where $\Natural^0 = \{\emptyset\}$, $\emptyset$ is the root of the tree and each vertex $\alpha=(n_1,\ldots,n_p)\in \Natural^{p}$ for $p\leq r-1$ has children 
$$
\alpha n : = (n_1,\ldots,n_p,n) \in \Natural^{p+1}
$$
for all $n\in \Natural$. Each vertex $\alpha$ is connected to the root $\emptyset$ by the path
$$
\emptyset \to n_1 \to (n_1,n_2) \to\cdots\to (n_1,\ldots,n_p) = \alpha.
$$
We will denote the set of vertices in this path by 
\begin{equation}
p(\alpha) = \bigl\{ \emptyset, n_1, (n_1,n_2),\ldots,(n_1,\ldots,n_p)  \bigr\}.
\label{pathtoleaf}
\end{equation}
We will consider rearrangements of $\Natural^r$ that preserve the structure of the tree $\A(r)$, in the sense that they preserve the parent-child relationship. More specifically, we define by
\begin{equation}
\alpha\wedge\beta 
:=
 |p(\alpha) \cap p(\beta)  |
\label{wedge}
\end{equation}
the number of common vertices in the paths from the root $\emptyset$ to the vertices $\alpha$ and  $\beta$, and consider the following group of maps on $\Natural^r$,
\begin{equation}
H_r = \bigl\{ \pi: \Natural^r\to \Natural^r \,\bigr|\, \pi \mbox{ is a bijection}, \pi(\alpha)\wedge \pi(\beta) = \alpha\wedge\beta \mbox{ for all } \alpha,\beta\in \Natural^r \bigr\}.
\label{setH}
\end{equation}
Any such map can be extended to the entire tree $\A(r)$ in a natural way: let $\pi(\emptyset) := \emptyset$ and
\begin{equation}
\mbox{if } \pi((n_1,\ldots,n_r)) = (m_1,\ldots,m_r) \,\mbox{ then let }\, \pi((n_1,\ldots,n_p)) := (m_1,\ldots,m_p).
\label{setHA}
\end{equation}
Because of the condition $\pi(\alpha)\wedge \pi(\beta) = \alpha\wedge\beta$ in (\ref{setH}), this definition does not depend on the coordinates $n_{p+1},\ldots,n_r$, so the extension is well-defined. It is clear that the extension preserves the parent-child relationship.  For each $\alpha \in \A(r)\setminus \Natural^r$, it follows that $\pi(\alpha n) = \pi(\alpha) \pi_\alpha(n)$ for some bijection $\pi_\alpha: \Natural \to \Natural$. In other words, the condition $\pi(\alpha)\wedge \pi(\beta) = \alpha\wedge\beta$ means that we can visualize the map $\pi$ as a recursive procedure, in which children $\alpha n$ of the vertex $\alpha\in \Natural^p$ are rearranged among themselves for each $\alpha$.  Note that $H_1$ is simply the group of all permutations of $\Natural$.

We will say that an array of random variables $(X_\alpha)_{\alpha\in \Natural^r}$ taking values in a standard Borel space $A$ (i.e. Borel-isomorphic to a Borel subset of a Polish space) is \emph{hierarchically exchangeable}, or \emph{$H$-exchangeable}, if 
\begin{equation}
\bigl(X_{\pi(\alpha)} \bigr)_{\alpha\in \Natural^r}
\stackrel{d}{=}
\bigl(X_\alpha \bigr)_{\alpha\in \Natural^r}
\label{HexchDF}
\end{equation}
for all $\pi\in H_r$. Throughout the paper, we will view any array of random variables as a random element in the product space, so the equality in distribution is always in the sense of equality of the finite dimensional distributions. Because of this, one can replace the condition in (\ref{setH}) that $\pi$ is a bijection by the condition that $\pi$ is simply an injection, since any injection viewed on finitely many elements can be, obviously, extended to a bijection preserving the property $\pi(\alpha)\wedge \pi(\beta) = \alpha\wedge\beta$.

The case of $r=1$ corresponds to the classical notion of an exchangeable sequence, and in the general case of $r\geq 1$ we will prove the following analogue of de Finetti's classical theorem. One natural example of an $H$-exchangeable array is given by (recall the notation in (\ref{pathtoleaf}))
\begin{equation}
X_\alpha = \sigma\bigl((v_{\beta})_{\beta\in p(\alpha)} \bigr),
\label{sigmaDF}
\end{equation}
where $\sigma: [0,1]^{r+1} \to A$ is a measurable function, and $v_\alpha$ for $\alpha\in\A(r)$ are i.i.d. random variables with the uniform distribution on $[0,1]$. The reason this array is hierarchically exchangeable is because, by the definition of $\pi$, the random variables $v_{\pi(\alpha)}$ for $\alpha\in \A(r)$ are also i.i.d. and uniform on $[0,1]$, $p(\pi(\alpha))= \pi(p(\alpha))$ and $X_{\pi(\alpha)} = \sigma((v_{\pi(\beta)})_{\beta\in p(\alpha)}).$ We will show the following.

\begin{theorem}\label{Th1}
Any hierarchically exchangeable array $(X_\alpha)_{\alpha\in \Natural^r}$ can be generated in distribution as in (\ref{sigmaDF}) for some measurable function $\sigma$.
\end{theorem}

This result is not very difficult to prove, and one can give several different arguments. We will describe an approach that will be a natural first step toward the general case of processes indexed by several trees or, more specifically, by product sets of the form $\Natural^{r_1}\times \cdots \times \Natural^{r_\ell}$ for any integers $r_1,\ldots,r_\ell \geq 1$. Recalling the definition (\ref{setH}), let us denote 
\begin{equation}
H_{r_1,\ldots,r_\ell} = H_{r_1}\times \cdots \times H_{r_\ell},
\label{Hmany}
\end{equation}
and for any $\pi=(\pi_1,\ldots,\pi_\ell)\in H_{r_1,\ldots,r_\ell}$ and any $\alpha = (\alpha_1,\ldots,\alpha_\ell) \in \Natural^{r_1}\times \cdots \times \Natural^{r_\ell}$, let us denote
$$
\pi(\alpha) = \bigl(\pi_1(\alpha_1),\ldots,\pi_\ell(\alpha_\ell) \bigr).
$$
We will say that an array of random variables $X_{\alpha}$ indexed by ${\alpha\in \Natural^{r_1}\times \cdots \times \Natural^{r_\ell}}$ and taking values in a standard Borel space $A$ is \emph{hierarchically exchangeable}, or \emph{$H$-exchangeable}, if
\begin{equation}
\bigl(X_{\pi(\alpha)} \bigr)_{\alpha\in \Natural^{r_1}\times \cdots \times \Natural^{r_\ell}}
\stackrel{d}{=}
\bigl(X_\alpha \bigr)_{\alpha\in \Natural^{r_1}\times \cdots \times \Natural^{r_\ell}}
\label{HexchDFmany}
\end{equation}
for all $\pi\in H_{r_1,\ldots,r_\ell}$. Let us denote
$$
\A(r_1,\ldots,r_\ell) = \A(r_1)\times \cdots \times \A(r_\ell)
$$
and, for $\alpha = (\alpha_1,\ldots,\alpha_\ell) \in\A(r_1,\ldots,r_\ell) $, denote
$$
p(\alpha) := p(\alpha_1)\times \cdots \times p(\alpha_\ell).
$$
Then, again, the natural class of $H$-exchangeable arrays is those of the form
\begin{equation}
X_\alpha = \sigma\bigl((v_{\beta})_{\beta\in p(\alpha)} \bigr),
\label{sigmaDFmany}
\end{equation}
for some measurable function $\sigma:[0,1]^{(r_1+1) + \ldots + (r_\ell+1)}\to A$ and a family of i.i.d. random variables $v_\beta$ indexed by ${\beta \in \A(r_1,\ldots,r_\ell)}$ with the uniform distribution on $[0,1]$.

\begin{theorem}\label{Th2}
Any hierarchically exchangeable array $(X_{\alpha})_{\alpha\in \Natural^{r_1}\times \cdots \times \Natural^{r_\ell}}$ can be generated in distribution as in (\ref{sigmaDFmany}) for some measurable function $\sigma$.
\end{theorem} 

This general result was motivated by the following special case, when the array $(X_{\alpha, i})$ is indexed by $\alpha\in \Natural^r$ and $i\in\Natural$. The condition (\ref{HexchDFmany}) now becomes 
\begin{equation}
\bigl(X_{\pi(\alpha), \rho(i)} \bigr)_{\alpha\in \Natural^r, i\in \Natural}
\stackrel{d}{=}
\bigl(X_{\alpha,i} \bigr)_{\alpha\in \Natural^r, i\in\Natural}
\label{HexchAH}
\end{equation}
for all $\pi\in H_r$ and all bijections $\rho:\Natural \to \Natural$, and Theorem \ref{Th2} implies that any such array can be generated in distribution as
\begin{equation}
X_{\alpha,i} = \sigma\bigl((v_{\beta})_{\beta\in p(\alpha)}, (v_{\beta}^i)_{\beta\in p(\alpha)} \bigr),
\label{sigmaAH}
\end{equation}
where $\sigma: [0,1]^{2(r+1)} \to\Reals$ is a measurable function and all $v_\alpha$ and $v_\alpha^i$ for $\alpha\in\A(r)$ and $i\in \Natural$ are i.i.d. random variables with the uniform distribution on $[0,1]$. This can be viewed as a hierarchical version of the classical \emph{Aldous-Hoover representation} (\cite{Aldous}, \cite{Aldous2}, \cite{Hoover}, \cite{Hoover2}), which corresponds to the case $r=1$. One application of this representation can be found in \cite{HEPS}, where it is explained how (\ref{sigmaAH}) is related to the predictions about the structure of the Gibbs measure in diluted spin glass models that originate in the work of M\'ezard and Parisi \cite{Mezard}. The main result in \cite{HEPS} proves precisely the hierarchical exchangeability (\ref{HexchAH}) for the random variables $X_{\alpha,i}$ that represent the magnetization of the $i^{\rm th}$ spin inside the pure state $\alpha$, and the tree structure as above stems from the ultrametric organization of the pure states in the Parisi ansatz, which was recently proved in \cite{PUltra}. Finally, although this is not directly related to the results presented in this paper, an interested reader can find a study of another notion of exchangeability on (infinite infinitary) trees in Section III.13 in \cite{Aldous2}. 

\medskip
\noindent
\textbf{Acknowledgement.} We would like to thank the referees for their careful review and a number of suggestions to improve the quality of the paper.

\section{The case of one tree}

It is well known that any standard Borel space is Borel-isomorphic to a Borel subset of $[0,1]$ (see e.g. Section $13.1$ in \cite{Dudley}), which means that it is enough to prove Theorems \ref{Th1} and \ref{Th2} with random variables $X_\alpha$ taking values in $[0,1]$, which we will assume from now on. All the arrays that we will deal with will take values in the product space of countably many copies of $[0,1]$, which is a compact space. For simplicity of notation, we will continue to denote all such spaces by $A$. We will denote by $\Pr\,A$ the space of probability measures on $A$ equipped with the topology of weak convergence, which is also a compact space. If a sequence $(X_n)_n$ of $A$-valued random variables is such that the empirical distributions
$$
\frac{1}{N}\sum_{n=1}^N\delta_{X_n}
$$
converge almost surely to some $(\Pr\,A)$-valued random variable, then we will call this limit the \emph{empirical measure} of $(X_n)_n$ and denote it by $\mathcal{E}((X_n)_n)$. Our key tool will be the following strong version of de Finetti's theorem (see Proposition 1.4, Corollary 1.5 and Corollary 1.6 from~\cite{Kallenberg}).

\begin{theorem}[de Finetti-Hewitt-Savage Theorem]\label{thm:deFHS}
Suppose $(X_n)_n$ is an exchangeable sequence of $A$-valued random variables. Then the empirical measure $\mathcal{E}((X_n)_n)$ exists almost surely and has the following properties:
\begin{itemize}
\item[(i)] $\mathcal{E}((X_n)_n)$ is almost surely a function of $(X_n)_n$;
\item[(ii)] given $\mathcal{E}((X_n)_n)$, the random variables $X_n$ are i.i.d. with the distribution $\mathcal{E}((X_n)_n)$;
\item[(iii)] if $Z$ is any other random variable on the same probability space such that 
\begin{eqnarray}\label{eq:cond-exch}
(Z,X_1,X_2,\ldots) \stackrel{d}{=} (Z,X_{\pi(1)},X_{\pi(2)},\ldots)
\,\,\mbox{ for all $\pi\in H_1$}
\end{eqnarray}
then the sequence $(X_n)_n$ is conditionally independent from $Z$ given $\mathcal{E}((X_n)_n)$. 
\end{itemize}
\end{theorem}

\medskip
\noindent
\textbf{Proof of Theorem \ref{Th1}.} The proof will be by induction on $r\geq 1$. For each $\alpha \in \Natural^{r-1}$, by Theorem~\ref{thm:deFHS}, the empirical measures
\begin{equation}
X_\alpha := \mathcal{E}\bigl((X_{\alpha n})_{n}\bigr) \in \Pr\,A
\end{equation}
exist almost surely, because hierarchical exchangeability (\ref{HexchDF}) implies that $(X_{\alpha n})_n$ is exchangeable in the index $n$ for each fixed $\alpha$. Moreover, hierarchical exchangeability together with Theorem~\ref{thm:deFHS} imply the following:
\begin{enumerate}
\item[(a)] Given $X_\alpha$ for a fixed $\alpha\in\Natural^{r-1}$, the random variables $X_{\alpha n}$, $n\in\Natural$, are i.i.d. with the distribution $X_\alpha$.
\item[(b)] The random variables $(X_{\alpha n})_{\alpha\in\Natural^{r-1}, n\in\Natural}$ are conditionally independent given $(X_\alpha)_{\alpha\in\Natural^{r-1}}$.  This holds because for a chosen $\alpha$, the joint distribution of all the random variables is invariant if one permutes the sequence $(X_{\alpha n})_{n \in \Natural}$ while leaving all $(X_{\alpha' n})_{\alpha'\neq \alpha,\,n\in\Natural}$ fixed, and so (iii) of Theorem 3 gives that the former are conditionally independent from the latter over $X_\alpha$.
\item[(c)] The empirical measures $(X_\alpha)_{\alpha\in\Natural^{r-1}}$ are hierarchically exchangeable,
$$
(X_{\pi(\alpha)})_{\alpha \in \Natural^{r-1}} \stackrel{d}{=} (X_\alpha)_{\alpha \in \Natural^{r-1}} \,\,\mbox{ for all }\, \pi \in H_{r-1}.
$$
\end{enumerate}
By the induction hypothesis, property (c) yields a representation
\begin{equation}
(X_\beta)_{\beta \in \Natural^{r-1}} \stackrel{d}{=} \big(\sigma_1((\nu_\gamma)_{\gamma \in p(\beta)})\big)_{\beta \in \Natural^{r-1}}.
\label{indhyp}
\end{equation}
By the properties (a) and (b) and the fact that $A$ is a Borel space, there exists a measurable function $\sigma_2:\Pr\,A\times [0,1]\to A$ such that, conditionally on $(X_\alpha)_{\alpha\in\Natural^{r-1}}$,
\begin{equation}
\bigl(X_{\alpha n}\bigr)_{\alpha\in\Natural^{r-1}, n\in\Natural} 
\stackrel{d}{=} \bigl(\sigma_2(X_\alpha,v_{\alpha n}) \bigr)_{\alpha\in\Natural^{r-1}, n\in\Natural},
\label{nolind}
\end{equation}
where $v_{\alpha n}$ for $\alpha n \in \Natural^r$ are i.i.d. random variables uniform on $[0,1]$, independent from everything else. In other words, we simply realize independent random variables $X_{\alpha n}$ from the distribution $X_\alpha$ as functions of independent uniform random variables $v_{\alpha n}$. (See, for instance, Lemma 7.8 in~\cite{Kallenberg} for a rather stronger result guaranteeing that this can be done.) Combining (\ref{indhyp}) and (\ref{nolind}) implies
\[(X_\alpha)_{\alpha \in \Natural^r} \stackrel{d}{=} \big(\sigma((\nu_\beta)_{\beta \in p(\alpha)})\big)_{\alpha \in \Natural^r}\]
with $\sigma(x_0,x_1,\ldots,x_r) := \sigma_2(\sigma_1(x_0,x_1,\ldots,x_{r-1}),x_r)$, which finishes the proof.
\qed

\section{The case of several trees}

Theorem~\ref{Th2} will be proved by induction on $(r_1,\ldots,r_\ell)$. Of course, the case $\ell = 1$ is already proved in the previous section.  However, in order to close the induction, it will actually be convenient to focus on a more general result, describing $H$-exchangeable couplings between processes and $I$-fields, defined as follows. We will call an array of random variables $(u_\alpha)_{\alpha \in \A(r_1,\ldots,r_\ell)}$ taking values in some compact spaces \emph{an $I$-field} if all $u_\alpha$ are independent and the distribution of $u_\alpha$ depends only on the ``distance of $\alpha$ from the root", namely, all $u_\alpha$ have the same distribution for $\alpha\in  \Natural^{p_1}\times \cdots \times \Natural^{p_\ell}$ for any given $(p_1,\ldots,p_\ell)$. We will consider a pair of processes
\begin{equation}
(u_\alpha)_{\alpha \in \A(r_1,\ldots,r_\ell)},(X_\alpha)_{\alpha\in \Natural^{r_1}\times \cdots \times \Natural^{r_\ell}}, 
\label{Upair}
\end{equation}
where $(u_\alpha)$ is an $I$-field, not necessarily independent of $(X_\alpha)$. We will assume that they are jointly hierarchically exchangeable in the sense that
\begin{equation}
\bigl((u_{\pi(\alpha)})_{\alpha \in \A(r_1,\ldots,r_\ell)}, (X_{\pi(\alpha)})_{\alpha\in \Natural^{r_1}\times \cdots \times \Natural^{r_\ell}}\bigr)
\stackrel{d}{=}
\bigl((u_{\alpha})_{\alpha \in \A(r_1,\ldots,r_\ell)}, (X_{\alpha})_{\alpha\in \Natural^{r_1}\times \cdots \times \Natural^{r_\ell}}\bigr)
\label{HexchUfield}
\end{equation}
for all bijections $\pi\in H_{r_1,\ldots, r_\ell}$ in (\ref{Hmany}) extended in a natural way to the entire set $ \A(r_1,\ldots,r_\ell)$, i.e. each coordinate $\pi_i\in H_{r_i}$ is extended from $\Natural^{r_i}$ to $\A(r_i)$ as in (\ref{setHA}). For convenience of notation, given an array $Y_\alpha$ indexed by $\alpha\in \A(r_1,\ldots,r_\ell)$ and a subset $S\subseteq \A(r_1,\ldots,r_\ell)$, we will denote $Y_S = (Y_\alpha)_{\alpha\in S}$. For example, $Y_{p(\alpha)} = (Y_\beta)_{\beta\in p(\alpha)}$. The following proposition is a generalization of Theorem \ref{Th2}.
\begin{proposition}\label{prop}
If (\ref{HexchUfield}) holds then there exists a measurable function $\tau$ such that, conditionally on the $I$-field $(u_\alpha)_{\alpha \in \A(r_1,\ldots,r_\ell)}$,
\begin{equation}
(X_\alpha )_{\alpha\in \Natural^{r_1}\times \cdots \times \Natural^{r_\ell}}
\stackrel{d}{=} 
\bigl(\tau(u_{p(\alpha)},v_{p(\alpha)}) \bigr)_{\alpha\in \Natural^{r_1}\times \cdots \times \Natural^{r_\ell}},
\label{propeq}
\end{equation}
where $(v_\alpha)_{\alpha \in \A(r_1,\ldots,r_\ell)}$ are i.i.d. random variables uniform on $[0,1]$, independent of $(u_\alpha)_\alpha$.
\end{proposition}
\medskip
\noindent
Formally, this equality of distribution conditionally on $(u_\alpha)_{\alpha \in \A(r_1,\ldots,r_\ell)}$ means the following equality of distribution for larger families of random variables:
\begin{eqnarray*}
&& \!\!\!
\Bigl((u_\alpha)_{\alpha \in \A(r_1,\ldots,r_\ell)},\ (X_\alpha )_{\alpha\in \Natural^{r_1}\times \cdots \times \Natural^{r_\ell}}\Bigr)
\\
&\stackrel{d}{=}& \!\!\!
\Bigl((u_\alpha)_{\alpha \in \A(r_1,\ldots,r_\ell)},\ \bigl(\tau(u_{p(\alpha)},v_{p(\alpha)}) \bigr)_{\alpha\in \Natural^{r_1}\times \cdots \times \Natural^{r_\ell}}\Bigr).
\end{eqnarray*}
We will generally avoid writing this out in full for the sake of lighter notation.

Of course, (\ref{propeq}) implies Theorem \ref{Th2} by considering an $I$-field $(u_\alpha)$ independent of the process $(X_\alpha)$. Proposition \ref{prop} will be proved by induction on $(r_1,\ldots,r_\ell)$ and, in the induction step, we will need to describe a conditional distribution of one array given another. We will be able to replace this second array with an $I$-field, and the independence built into the definition of $I$-fields will be well-suited for the induction argument. The induction argument does not work so well when the $I$-field in Proposition \ref{prop} is replaced by a general $H$-exchangeable array $(Y_\alpha)$. However, such a generalization, described in Theorem \ref{Th4} below, will follow once we have Proposition \ref{prop}. 

To describe the induction, it will be convenient to write members of $\A(r_1,\ldots,r_\ell)$ in the form $(\omega,\alpha)$, where $\omega \in \A(r_1,\ldots,r_{\ell-1})$ and $\alpha \in \A(r_\ell)$, and also abbreviate
$$
\A = \A(r_1,\ldots,r_{\ell-1})
\,\mbox{ and }\,
\LL = \Natural^{r_1}\times \cdots \times \Natural^{r_{\ell-1}}.
$$
We therefore write the pair of processes (\ref{Upair}) as 
$
(u_{\omega,\alpha})_{\omega\in \A,\alpha \in \A(r_\ell)},(X_{\omega,\alpha})_{\omega \in \LL,\alpha \in \Natural^{r_\ell}}.
$
To close the induction we will make three separate appeals to simpler cases of Proposition~\ref{prop}, and we subdivide the proof into stages accordingly.

\subsubsection*{Using the case of one tree}

\noindent
For the first stage, it will also be convenient to introduce the notation, for each $\alpha \in \Natural^{r_\ell}$,
\begin{equation}
\t{X}_\alpha = (\t{X}_\alpha^1, \t{X}_\alpha^2)
=
\bigl((u_{\omega,\alpha})_{\omega \in \A},(X_{\omega,\alpha})_{\omega \in \LL}\bigr),
\label{tildeX}
\end{equation}
which is an element of another compact space, say $\t{A} = \t{A}_1 \times \t{A}_2$, where $\t{X}_\alpha^j$ take values in $\t{A}_j$ for $j=1,2.$ If we denote the subarray 
\begin{equation}
U^- = (u_{\omega,\alpha})_{\omega\in \A, \alpha\in \A(r_\ell-1)}
\label{Uminus}
\end{equation}
of our $I$-field consisting of the coordinates that do not appear in (\ref{tildeX}), then in these terms our goal is to describe the joint distribution of $(\t{X}_\alpha)_{\alpha\in \Natural^{r_\ell}}$ and $U^-$.

First of all, notice that hierarchical exchangeability in (\ref{HexchUfield}) implies that the process $(\t{X}_\alpha)_{\alpha\in \Natural^{r_\ell}}$ is $H$-exchangeable. Hence, similarly to the proof of Theorem \ref{Th1}, for each $\alpha \in \Natural^{r_\ell-1}$, the empirical measure
\begin{equation}
\t{X}_\alpha := \mathcal{E}\bigl((\t{X}_{\alpha n})_{n} \bigr) \in \Pr\,\t{A}
\label{tildeEmp}
\end{equation}
exists almost surely and, by Theorem~\ref{thm:deFHS}, we get:
\begin{enumerate}
\item[(a)] given $\t{X}_\alpha$ for $\alpha\in\Natural^{r_\ell-1}$, the random variables $\t{X}_{\alpha n}$ are i.i.d. with the distribution $\t{X}_\alpha$;
\item[(b)] given $(\t{X}_\alpha)_{\alpha\in\Natural^{r_\ell-1}}$, the random variables $(\t{X}_{\alpha n})_{\alpha\in\Natural^{r_\ell-1}, n\in\Natural}$ are conditionally independent.
\end{enumerate}
Note also that the permutation of the index $n$ for a fixed $\alpha$ does not affect the subarray (\ref{Uminus}). Therefore, part (iii) of Theorem~\ref{thm:deFHS} also implies that 
\begin{enumerate}
\item[(c)] given $(\t{X}_\alpha)_{\alpha\in\Natural^{r_\ell-1}}$, the array $(\t{X}_{\alpha n})_{\alpha\in\Natural^{r_\ell-1}, n\in\Natural}$ is independent of $U^-$. 
\end{enumerate}
Another important observation is that, by the definition of $I$-field, for any $\alpha\in\Natural^{r_\ell-1}$, the random variables $\t{X}_{\alpha n}^1 = (u_{\omega,\alpha n})_{\omega \in \A}$ in (\ref{tildeX}) are i.i.d. for $n\in \Natural$ with some fixed distribution on $\t{A}_1$ and, therefore, the marginal of the empirical measure $\t{X}_\alpha$ in (\ref{tildeEmp}) on $\t{A}_1$ is this fixed nonrandom measure. Together with the property (a) this implies:
\begin{enumerate}
\item[(d)]
the random variables $\t{X}_{\alpha n}^1$ for $n\in \Natural$ are independent of the empirical measure $\t{X}_\alpha$.
\end{enumerate}
Let us now consider an infinite subset $I\subseteq \Natural$ such that $I^c=\Natural\setminus I$ is also infinite.  Even though our goal is to describe the joint distribution of $(\t{X}_{\alpha})_{\alpha\in \Natural^{r_\ell}}$ and $U^-$, because of the hierarchical exchangeability it is, obviously, sufficient to describe the joint distribution of 
$$
(\t{X}_{\alpha n})_{\alpha\in \Natural^{r_\ell-1}, n\in I} \,\,\mbox{ and }\,\, U^-.
$$ 
This will be done in several steps, and we begin with the following lemma. We will suppose, without loss of generality,  that $1\in I$. We will write $\p(Y\in \cdot\,|\,Y')$ for the conditional distribution of $Y$ given $Y'$.
\begin{lemma} \label{LemABC}(A) The following equality holds:
\begin{eqnarray}
\p\Bigl( (\t{X}_{\alpha n})_{\alpha\in\Natural^{r_\ell-1}, n\in I} \in\ \cdot\ \Big|\ 
(\t{X}_{\alpha n})_{\alpha\in\Natural^{r_\ell-1}, n\in I^c}\Bigr)
=\,\,
\bigotimes\nolimits_{\alpha\in\Natural^{r_\ell-1}}
\p\Bigl( \t{X}_{\alpha 1} \in\ \cdot\ \Big|\ (\t{X}_{\alpha n})_{ n\in I^c}\Bigr)^{\otimes I}.
\label{abc}
\end{eqnarray}
(B) Conditionally on $(\t{X}_{\alpha n})_{\alpha\in\Natural^{r_\ell-1}, n\in I^c}$, the arrays $(\t{X}_{\alpha n})_{\alpha\in\Natural^{r_\ell-1}, n\in I}$ and $U^-$ are independent.

\noindent
(C) The arrays $(\t{X}_{\alpha n}^1)_{\alpha\in\Natural^{r_\ell-1}, n\in I}$ and $(\t{X}_{\alpha n})_{\alpha\in\Natural^{r_\ell-1}, n\in I^c}$ are independent.
\end{lemma}
\textbf{Proof.} First of all, by property (a), the empirical measure (\ref{tildeEmp}) satisfies
\begin{equation}
\t{X}_\alpha = \mathcal{E}\bigl((\t{X}_{\alpha n})_{n\in I^c} \bigr), 
\label{tildeEmpI}
\end{equation}
which means that $\t{X}_\alpha$ is almost surely a function of $(\t{X}_{\alpha n})_{n\in I^c}$. Therefore,
\begin{eqnarray*}
&&
\p\Bigl( (\t{X}_{\alpha n})_{\alpha\in\Natural^{r_\ell-1}, n\in I} \in\ \cdot\ \Big|\ 
(\t{X}_{\alpha n})_{\alpha\in\Natural^{r_\ell-1}, n\in I^c},\ U^-\Bigr)
\\
&& =\,\,
\p\Bigl( (\t{X}_{\alpha n})_{\alpha\in\Natural^{r_\ell-1}, n\in I} \in\ \cdot\ \Big|\ 
(\t{X}_{\alpha})_{\alpha\in\Natural^{r_\ell-1}},
(\t{X}_{\alpha n})_{\alpha\in\Natural^{r_\ell-1}, n\in I^c},\ U^-\Bigr).
\end{eqnarray*} 
Using the properties (b) and (c), this conditional distribution is equal to
\begin{equation}
\p\Bigl( (\t{X}_{\alpha n})_{\alpha\in\Natural^{r_\ell-1}, n\in I} \in\ \cdot\ \Big|\ 
(\t{X}_\alpha)_{\alpha\in\Natural^{r_\ell-1}}\Bigr).
\label{abcinter}
\end{equation}
The same computation obviously also works without $U^-$, and therefore
\begin{eqnarray*}
&&
\p\Bigl( (\t{X}_{\alpha n})_{\alpha\in\Natural^{r_\ell-1}, n\in I} \in\ \cdot\ \Big|\ 
(\t{X}_{\alpha n})_{\alpha\in\Natural^{r_\ell-1}, n\in I^c},\ U^-\Bigr)
\\
&& =\,\,
\p\Bigl( (\t{X}_{\alpha n})_{\alpha\in\Natural^{r_\ell-1}, n\in I} \in\ \cdot\ \Big|\ 
(\t{X}_{\alpha n})_{\alpha\in\Natural^{r_\ell-1}, n\in I^c}\Bigr).
\end{eqnarray*} 
This proves (B). Next, using the properties (a) and (b), we can rewrite (\ref{abcinter}) as (recall that $1\in I$)
$$
\bigotimes\nolimits_{\alpha\in\Natural^{r_\ell-1}, n\in I} \,\,
\p\Bigl(\t{X}_{\alpha n}\in\ \cdot\ \Big|\ \t{X}_\alpha \Bigr)
=
\bigotimes\nolimits_{\alpha\in\Natural^{r_\ell-1}}
\p\Bigl(\t{X}_{\alpha 1}\in\ \cdot\ \Big|\ \t{X}_\alpha \Bigr)^{\otimes I},
$$
which proves that
\begin{equation}
\p\Bigl( (\t{X}_{\alpha n})_{\alpha\in\Natural^{r_\ell-1}, n\in I} \in\ \cdot\ \Big|\ 
(\t{X}_{\alpha n})_{\alpha\in\Natural^{r_\ell-1}, n\in I^c}\Bigr)
=
\bigotimes\nolimits_{\alpha\in\Natural^{r_\ell-1}}
\p\Bigl(\t{X}_{\alpha 1}\in\ \cdot\ \Big|\ \t{X}_\alpha \Bigr)^{\otimes I}.
\label{abcproof}
\end{equation} 
Using (\ref{tildeEmpI}) and property (a), for any fixed $\alpha\in\Natural^{r_\ell-1}$,
\begin{eqnarray*}
&&
\p\Bigl( \t{X}_{\alpha 1} \in\ \cdot\ \Big|\ (\t{X}_{\alpha n})_{ n\in I^c}\Bigr)
=
\p\Bigl( \t{X}_{\alpha 1} \in\ \cdot\  \Big|\  \t{X}_\alpha, (\t{X}_{\alpha n})_{ n\in I^c}\Bigr)
=
\p\Bigl( \t{X}_{\alpha 1} \in\ \cdot\ \Big|\ \t{X}_\alpha\Bigr).
\end{eqnarray*} 
Combining the last two equations proves (A). The last claim follows from (\ref{abcproof}) and property (d) above.
\qed

\subsubsection*{Using the case of $\A\times \A(r_\ell-1)$}

\noindent
Now that we have utilized the exchangeability with respect to the permutations of the index $n$, we will change the focus and make the dependence of all random variables on the index $\omega\in\A$ explicit. For each $\alpha\in \Natural^{r_\ell-1}$, let us denote
\begin{eqnarray}
X^{+}_{\omega,\alpha} \!\!&:=&\!\! (X_{\omega,\alpha n})_{n\in I^c} \,\mbox{ for }\, \omega\in \LL,
\label{switch1}
\\
U^+_{\omega,\alpha} \!\!&:=&\!\! (u_{\omega,\alpha n})_{n\in I^c} \,\mbox{ for }\, \omega\in \A,
\label{switch2}
\\
(U^+_\alpha, X^{+}_\alpha) \!\!&:=&\!\! \bigl((U^+_{\omega,\alpha})_{\omega \in \A},(X^{+}_{\omega,\alpha})_{\omega \in \LL}\bigr),
\label{switch3}
\\
(u_{\alpha n},X_{\alpha n}) \!\!&:=&\!\! \bigl((u_{\omega,\alpha n})_{\omega \in \A},(X_{\omega,\alpha n})_{\omega \in \LL}\bigr)
\,\mbox{ for }\, n\in I,
\label{switch4}
\end{eqnarray}
and let us also denote
\begin{eqnarray}
(U^+, X^{+}) \!\!&:=&\!\! (U^+_\alpha, X^{+}_\alpha)_{\alpha\in \Natural^{r_\ell-1}},
\label{switch5}
\\
(u,X) \!\!&:=&\!\! \bigl((u_{\alpha n},X_{\alpha n})\bigr)_{\alpha\in\Natural^{r_\ell-1}, n\in I}.
\label{switch}
\end{eqnarray}
With this notation, we can rewrite (\ref{abc}) as 
\begin{equation}
\p\Bigl( (u,X) \in\ \cdot\ \Big|\ 
(U^+, X^{+}) \Bigr)
=
\bigotimes\nolimits_{\alpha\in\Natural^{r_\ell-1}}
\p\Bigl(  (u_{\alpha 1},X_{\alpha 1}) \in\ \cdot\ \Big|\ 
(U^+_\alpha, X^{+}_\alpha) \Bigr)^{\otimes I}.
\label{abcswitch}
\end{equation}
We can also rewrite claims (B) and (C) in Lemma \ref{LemABC} as follows:
\begin{enumerate}
\item[(B$^\prime$)] conditionally on $(U^+, X^{+})$ the arrays $(u,X)$ and $U^-$ are independent;
\item[(C$^\prime$)] The arrays $u$ and $(U^+, X^{+})$ are independent.
\end{enumerate}
We will now make our first appeal to the inductive hypothesis of Proposition~\ref{prop} to describe the joint distribution of $(U^+, X^{+})$ and $U^-$. Notice that $U_{\omega,\alpha}^+$ in (\ref{switch2}) and some of the coordinates $u_{\omega,\alpha}$ in $U^-$ in (\ref{Uminus}) are indexed by $\omega\in \A$ and $\alpha\in \Natural^{r_\ell-1}$, so we will combine them and introduce a new array $U=(U_{\omega, \alpha})_{\omega\in \A, \alpha\in \A(r_\ell-1)}$ such that
\begin{eqnarray}
U_{\omega,\alpha} \!\! &: =&\!\! (u_{\omega,\alpha}, U_{\omega,\alpha}^+)
\,\mbox{ for }\,
\omega\in \A, \alpha\in \Natural^{r_\ell-1},
\label{UUplus}
\\
U_{\omega,\alpha} \!\! &: =&\!\! u_{\omega,\alpha}
\,\mbox{ for }\,
\omega\in \A, \alpha\in \A(r_\ell-1)\setminus \Natural^{r_\ell-1}.
\nonumber
\end{eqnarray}
Slightly abusing notation, this definition can be written as $U = (U^-, U^+)$ and it is obvious that $U$ is again an $I$-field. Let us also observe right away that, by property (B$^\prime$), 
\begin{equation}
\p\Bigl( (u,X) \in\ \cdot\ \Big|\ 
(U^+, X^{+}) \Bigr)
=
\p\Bigl( (u,X) \in\ \cdot\ \Big|\ 
(U, X^{+}) \Bigr).
\label{abcswitchU}
\end{equation}
The following gives a description of the joint distribution of $(U^+, X^{+})$ and $U^-$.
\begin{lemma}
Conditionally on the $I$-field $U= (U^-, U^+)$,
\begin{equation}
X^{+} =
\bigl(X^{+}_{\omega,\alpha}\bigr)_{\omega\in \LL,\alpha\in \Natural^{r_\ell-1}} 
\stackrel{d}{=}
\Bigl(\xi\bigl(v_{p(\omega,\alpha)}, U_{p(\omega,\alpha)} \bigr) \Bigr)_{\omega\in \LL, \alpha\in \Natural^{r_\ell-1}}
\label{induction1}
\end{equation}
for some measurable function $\xi$ of its coordinates, where $v_\beta$ are i.i.d. uniform random variables on $[0,1]$ indexed by $\A\times \A(r_\ell-1)$.
\end{lemma}
\textbf{Proof.} This is a consequence of the fact that $U$ is an $I$-field, and the pair $U$ and $(X^{+}_{\omega,\alpha})_{\omega\in \LL,\alpha\in \Natural^{r_\ell-1}}$ is, clearly, a hierarchically exchangeable coupling satisfying (\ref{HexchUfield}) with $r_\ell$ replaced by $r_{\ell}-1$. By the induction hypothesis, the claim follows.
 \qed

\medskip
\medskip
\noindent
Let us denote the array of random variables $v$ on the right hand side of (\ref{induction1}) by
$$
V
:=
\bigl(v_{\omega,\alpha}\bigr)_{\omega\in \A, \alpha\in \A(r_\ell-1)}. 
$$
Let us denote by ${\Xi}$ the full map on the right hand side of the equation (\ref{induction1}), which can be then written as 
$$
X^{+}\stackrel{d}{=}{\Xi}(V,U).
$$
Since all our random variables take values in standard Borel (or even compact) spaces, we can consider the regular conditional probability 
\begin{equation}
\p\bigl(\ \cdot \ \bigr| \ x\, \bigr)
=
\p\Bigl(
V \in\ \cdot\ \Big|\  \bigl(U, \Xi(V,U)\bigr) = x
\Bigr).
\label{condprob}
\end{equation}
It is a standard fact in this case that if $\mu$ is the law of $(U, X^{+})$ then, for $\mu$-almost all $x$,
\begin{equation}
\p\Bigl(
\bigl\{V \, \bigr| \, (U, {\Xi}(V,U)) = x \bigr\} \ \Big|\  x
\Bigr)=1.
\label{condprob2}
\end{equation}
Now, using this conditional probability, let us couple the arrays $(u,X)$ in (\ref{switch}) and $V$ conditionally independently given $(U,X^{+})$,
\begin{eqnarray}
&&
\p\Bigl(
(u,X), V 
 \in\ \cdot\ \Big|\ 
(U,X^{+})=x
\Bigr)
\nonumber
\\
&&
=\,\,
\p\Bigl(
(u,X)
 \in\ \cdot\ \Big|\ 
(U,X^{+})=x
\Bigr)
\times
\p\Bigl(
V
 \in\ \cdot\ \Big|\ 
(U,X^{+})=x
\Bigr).
\label{join}
\end{eqnarray}
This is a standard construction in probability, as well as in ergodic theory, where it is called a `relatively independent joining': see, for instance, the third example in Section 6.1 of Glasner~\cite{Gla}.  The triple
$$
(u,X), V \mbox{ and } (U,X^{+})
$$
is still hierarchically exchangeable, since this is true separately of both conditional distributions on the right hand side of (\ref{join}) (for a much more detailed explanation see Lemma 2.3 in \cite{Kallenberg0}).  Having done this, we may henceforth regard all of these processes as defined on the same background probability space.

\begin{lemma}\label{LemInd2}
With the joint distribution constructed above,
\begin{equation}
\p\Bigl((u,X)\in\ \cdot\ \Big|\ (U^+,X^{+}) \Bigr)
= \p\Bigl((u,X)\in\ \cdot\ \Big|\ (V,U)\Bigr).
\label{xitoVU}
\end{equation}
\end{lemma}
Notice that this implies that the property (C$^\prime$) above can now be written as: 
\begin{enumerate}
\item[(C$^\prime$$^\prime$)] the arrays $u$ and $(V, U)$ are independent.
\end{enumerate}
\medskip
\noindent
\textbf{Proof of Lemma \ref{LemInd2}.} By (\ref{condprob2}), $X^{+} = {\Xi}(V,U)$ with probability one, so $X^{+}$ is almost surely a function of $V$ and $U$. Therefore,
$$
\p\Bigl((u,X)\in\ \cdot\ \Big|\ (V,U)\Bigr) = \p\Bigl((u,X)\in\ \cdot\ \Big|\ X^{+}, (V,U)\Bigr).
$$ 
By the construction (\ref{join}), $(u,X)$ and $V$ are conditionally independently given $(U,X^{+})$, so this conditional distribution is equal to $\p\bigl((u,X)\in\ \cdot\ \big|\ (U,X^{+}) \bigr)$, and (\ref{abcswitchU}) finishes the proof.
\qed

\medskip
\medskip
\noindent
Thus, we have replaced the conditioning on $(U^+,X^{+})$ on the left hand side of (\ref{abcswitch}) with conditioning on $(V, U)$, and now we will do a similar substitution in each factor on the right hand side of (\ref{abcswitch}). Recall the notation $U_\alpha^+$ and $X_\alpha^+$ in (\ref{switch3}) and, for each $\alpha\in \Natural^{r_\ell-1}$, let us denote 
\begin{equation}
V_\alpha \!\!: =\!\! \bigl(v_{p(\omega,\alpha)}\bigr)_{\omega\in \LL}
\,\,\mbox{ and }\,\,
U_\alpha \!\!:=\!\! \bigl(U_{p(\omega,\alpha)} \bigr)_{\omega\in \LL}.
\label{switch8}
\end{equation}
Notice that one factor on the right hand side of (\ref{abcswitch}) is $\p\bigl((u_{\alpha 1}, X_{\alpha 1})\in\ \cdot \ \bigr| \ (U^+_\alpha,X^{+}_\alpha) \bigr)$ and we will now show the following.
\begin{lemma}
For each $\alpha\in \Natural^{r_\ell-1}$, we have
\begin{equation}
\p\Bigl((u_{\alpha 1}, X_{\alpha 1})\in\ \cdot \ \Bigr| \ (U^+_\alpha, X^{+}_\alpha) \Bigr)
=
\p\Bigl((u_{\alpha 1}, X_{\alpha 1})\in\ \cdot \ \Bigr| \ (V_\alpha, U_\alpha) \Bigr).
\end{equation}
\end{lemma}
\textbf{Proof.} First of all, the equation (\ref{abcswitch}) implies that
$$
\p\Bigl((u_{\alpha 1}, X_{\alpha 1})\in\ \cdot \ \Bigr| \ (U^+_\alpha, X^{+}_\alpha) \Bigr)
=
\p\Bigl((u_{\alpha 1}, X_{\alpha 1})\in\ \cdot \ \Bigr| \ (U^+, X^{+}) \Bigr),
$$
which can be seen by considering the probabilities of cylindrical sets that depend only on $(u_{\alpha 1}, X_{\alpha 1})$. Using (\ref{xitoVU}), we get
\begin{equation}
\p\Bigl((u_{\alpha 1}, X_{\alpha 1})\in\ \cdot \ \Bigr| \ (U^+_\alpha,X^{+}_\alpha) \Bigr)
=
\p\Bigl((u_{\alpha 1}, X_{\alpha 1})\in\ \cdot \ \Bigr| \ (V, U) \Bigr).
\label{extreme}
\end{equation}
We saw in the proof of Lemma \ref{LemInd2} that $X^{+} = {\Xi}(V,U)$ with probability one and, therefore,
$$
X^{+}_\alpha = \bigl(X^+_{\omega,\alpha}\bigr)_{\omega\in \LL} 
= \Bigl(\xi\bigl(v_{p(\omega,\alpha)}, U_{p(\omega,\alpha)} \bigr) \Bigr)_{\omega\in \LL}.
$$ 
Using this and the fact that, by (\ref{UUplus}), $U_\alpha^+$ is also a function of $U_\alpha$, we obtain the following inclusion of $\sigma$-algebras,
$$
\sigma(U_\alpha^+, X^{+}_\alpha)\subseteq \sigma(V_\alpha, U_\alpha) \subseteq \sigma(V, U).
$$
The equality of conditional distributions in (\ref{extreme}) given the two extreme $\sigma$-algebras implies the equality to the  conditional distribution given the middle $\sigma$-algebra, and this finishes the proof.
\qed

\medskip
\medskip
\noindent
The preceding two lemmas allow us to rewrite (\ref{abcswitch}) as
\begin{equation}
\p\Bigl( (u, X)\in\ \cdot \ \Bigr| \ (V, U)\Bigr)
=
\bigotimes\nolimits_{\alpha\in\Natural^{r_\ell-1}}
\p\Bigl( (u_{\alpha 1}, X_{\alpha 1})\in\ \cdot \ \Bigr| \ (V_\alpha, U_\alpha) \Bigr)^{\otimes I}.
\label{abcswitchVU}
\end{equation}
In other words, conditionally on $(V, U)$, the random variables $ (u_{\alpha n}, X_{\alpha n})$ are independent for all $\alpha\in\Natural^{r_\ell-1}$ and $n\in I$, and for a fixed $\alpha$, have the same distribution,
$$
\p\Bigl( (u_{\alpha 1}, X_{\alpha 1})\in\ \cdot \ \Bigr| \ (V_\alpha, U_\alpha) \Bigr).
$$
By the property (C$^\prime$$^\prime$) above, $u_{\alpha 1}$ is independent of $(V_\alpha, U_\alpha)$, so our main concern now is to describe the conditional distribution of $X_{\alpha 1}$ given $u_{\alpha 1},$ $V_\alpha$ and $U_\alpha$.

\subsubsection*{Using the case of $\ell-1$ trees}

\noindent
Lastly, we will use the induction hypothesis in Proposition \ref{prop} to describe the joint distribution of the processes $X_{\alpha 1}, u_{\alpha 1}, V_\alpha$ and $U_\alpha$  for a fixed  $\alpha\in \Natural^{r_\ell-1}$, so these are indexed by $\A = \A(r_1,\ldots,r_{\ell-1})$. The process $X_{\alpha 1}$ consists of the random variables $X_{\omega,\alpha 1}$ indexed by $\omega\in \LL.$ We will view the triple $(u_{\alpha 1}, V_\alpha,U_\alpha)$ as a new $I$-field that consists of the random variables
\begin{equation}
T_\omega^\alpha := \Bigl(u_{\omega, \alpha 1}, \bigl(v_{(\omega,\beta)}\bigr)_{\beta\in p(\alpha)}, 
\bigl(U_{(\omega,\beta)} \bigr)_{\beta\in p(\alpha)} \Bigr)
\end{equation}
indexed by $\omega\in \A$. Here, we relabeled the random variables by collecting all the coordinates of $U_\alpha$ and $V_\alpha$ that depend on a fixed $\omega\in \A.$ By the property (C$^\prime$$^\prime$) above, the array $T^\alpha:= (T^\alpha_\omega)_{\omega\in \A}$ is again an $I$-field, and it is clear that it forms a hierarchically exchangeable coupling with the array $X_{\alpha 1}$. The induction hypothesis in Proposition \ref{prop}, now used with $r_\ell=0$, implies the following.

\begin{lemma}\label{LemZeta}
There exists a measurable function $\tau$ such that, conditionally on $T^\alpha$,
\begin{equation}
\big(X_{\omega,\alpha 1}\bigr)_{\omega\in \LL}
\stackrel{d}{=} 
\bigl(\tau(w_{p(\omega)},T^\alpha_{p(\omega)})\bigr)_{\omega\in \LL},
\label{zetaw}
\end{equation}
where $w$ is an array of i.i.d. random variables uniform on $[0,1]$ indexed by $\omega\in \A$, independent of everything else.
\end{lemma}
This allows us to finish the proof of Proposition \ref{prop}. First of all, let us notice that we can write
$$
T^\alpha_{p(\omega)} = \bigl(u_{p(\omega)\times\{\alpha 1\}}, v_{p(\omega,\alpha)}, 
U_{p(\omega,\alpha)} \bigr).
$$
Combining Lemma \ref{LemZeta} with (\ref{abcswitchVU}), we proved that, conditionally on the arrays $u, V$ and $U$, we can generate the random variables $X_{\omega, \alpha n}$ for $\omega\in \LL, \alpha\in \Natural^{r_\ell-1}, n\in I$ in distribution by
\begin{equation}
X_{\omega, \alpha n} =
\tau\bigl(v_{p(\omega)\times \{\alpha n\}},u_{p(\omega)\times\{\alpha n\}}, v_{p(\omega,\alpha)}, 
U_{p(\omega,\alpha)} \bigr),
\label{tauend}
\end{equation}
where, for each $\alpha\in \Natural^{r_\ell-1}$ and $n\in I$, we used the random variables $v_{p(\omega)\times \{\alpha n\}}$ in place of an independent copy of $w_{p(\omega)}$ in (\ref{zetaw}). First of all, 
$$
\bigl(v_{p(\omega)\times \{\alpha n\}},  v_{p(\omega,\alpha)} \bigr) = v_{p(\omega, \alpha n)}.
$$
If we recall the definition of the process $U$ in (\ref{UUplus}), we see that for $\alpha\in\Natural^{r_\ell-1}$, $U_{p(\omega,\alpha)}$ consists of two parts, $u_{p(\omega,\alpha)}$ and $U^+_{p(\omega)\times\{\alpha\}}$, and the first one can be combined with $u_{p(\omega)\times \{\alpha n\}}$ to give
$$
\bigl(u_{p(\omega)\times \{\alpha n\}},  u_{p(\omega,\alpha)} \bigr) = u_{p(\omega,\alpha n)}.
$$
Then, (\ref{tauend}) can be rewritten as (slightly abusing notation)
\begin{equation}
X_{\omega, \alpha n} =
\tau\bigl(u_{p(\omega, \alpha n)}, v_{p(\omega, \alpha n)}, U^+_{p(\omega)\times\{\alpha\}} \bigr).
\label{tauend2}
\end{equation}
Finally, note that we consider the random variables $X_{\omega, \alpha n}$ with the index $n\in I$, while all the random variables $U^+_{\omega, \alpha}$ in (\ref{switch2}) were defined in terms of the random variables $u_{\omega,\alpha n}$ with the index $n\in I^c$, so now they are not viewed as a part of our $I$-field $(u,U^-)$. Therefore, by  redefining the function $\tau$, we can absorb the randomness of $U^+_{p(\omega)\times\{\alpha\}}$ into $v_{p(\omega, \alpha n)}$ to get
\begin{equation}
X_{\omega, \alpha n} =
\tau\bigl(u_{p(\omega, \alpha n)}, v_{p(\omega, \alpha n)} \bigr).
\label{tauend3}
\end{equation}
This completes the induction step in Proposition \ref{prop}, and finishes the proof of Theorem \ref{Th2}.
\qed

\medskip
\medskip
\noindent
One can also now formulate a conditional version of Theorem \ref{Th2} as follows. Examples of $H$-exchangeable pairs of processes can be constructed in the form
\begin{equation}
(Y_\alpha,X_\alpha) = \bigl(\sigma_1(u_{p(\alpha)}),\sigma_2(u_{p(\alpha)},v_{p(\alpha)})\bigr),
\label{pairsigmaAH}
\end{equation}
for two measurable functions $\sigma_1, \sigma_2$ and independent $I$-fields $u$ and $v$ of uniform random variables on $[0,1]$.

\begin{theorem}\label{Th4}
Any hierarchically exchangeable array of pairs $(Y_\alpha,X_\alpha)_{\alpha\in \Natural^{r_1}\times \cdots \times \Natural^{r_\ell}}$ can be generated in distribution as in (\ref{pairsigmaAH}) for some measurable functions $\sigma_1$ and $\sigma_2$.
\end{theorem} 
\textbf{Proof.} This follows by first applying Theorem~\ref{Th2} to represent
$$
(Y_\alpha)_{\alpha} \stackrel{d}{=} \bigl(\sigma_1(u_{p(\alpha)}) \bigr)_\alpha,
$$
then forming the coupling of the processes $X$ and $u$ conditionally independently over $Y$, and then applying Proposition~\ref{prop} to represent the joint distribution of $(u,X)$. \qed

\end{document}